\documentclass[11pt]{amsart}
\usepackage{amsthm,amssymb}
 
\newtheorem{thm}{Theorem}[section]
\newtheorem{lm}[thm]{Lemma}

\input epsf

\def \Z {\mathbb{Z}}

\def \M {\mathcal M}

\def \H {\mathcal H}
\def \P {\mathcal P}

\begin{document}

\title[]{Finite Type Link-homotopy Invariants}
\author[]{Xiao-Song Lin}
\address{Department of Mathematics, University of California, 
Riverside, CA 92521}
\email{xl@math.ucr.edu}
\thanks{Partially supported by the Overseas Youth Cooperation Research
Fund of NSFC}
\begin{abstract}{An explicit polynomial in the linking numbers
$l_{ij}$ and Milnor's triple linking
numbers $\mu(rst)$ on six component links is shown to be a well-defined
finite type link-homotopy invariant. This solves a problem 
raised by 
B. Mellor and D. Thurston. An extension of our construction
also produces a finite type
link invariant which detects the invertibility for some links.}  
\end{abstract}
\maketitle

\section{Introduction}

The classification of links in the 3-space up to link-homotopy \cite{HL} 
was published ten years ago. Since then, the question of whether one could 
extract link-homotopy invariants from this classification has not 
been addressed properly. Recall that this classification starts with
the classification of $k$ component string links up to link-homotopy by a 
finitely
generated torsion free nilpotent group $\H(k)$. Then link-homotopy classes are
classified as orbits of this group $\H(k)$ under the \lq\lq nilpotent action" of
conjugations and partial conjugations. The group $\H(k)$ is of rank
$$\sum_{n=2}^k(n-2)!\binom{k}{n}$$
so an element of $\H(k)$ can be described uniquely by that many integers. These
integers are Milnor's $\mu$-numbers\footnote{Usually, they are called $\mu$-invariants.
But the word \lq\lq invariant" is clearly abused here so we decide to call
them $\mu$-numbers.} with distinct indices. By a 
{\it link-homotopy invariant polynomial}, or simply a link-homotopy invariant, 
we mean a polynomial in these $\mu$-numbers 
which is invariant under the action of conjugations and partial conjugations.
There are trivial examples of such link-homotopy invariant polynomials coming
from polynomials of linking numbers. A link-homotopy invariant
polynomial is non-trivial if it contains higher order $\mu$-numbers.

The main result of this paper is that such a non-trivial 
link-homotopy invariant polynomial exists when 
$k\geq6$. 

The abelianization of $\H(k)$ is a free abelian group of rank $\binom{k}{2}$. This
is where the classical linking numbers $l_{ij}$, $1\leq i<j\leq k$, fit in. 
The action of conjugations and partial conjugations on this quotient of $\H(k)$
is trivial.
The next successive quotient of the lower central series of $\H(k)$ is a
free abelian group of rank $\binom{k}{3}$, whose elements can be described by
the collection of Milnor's triple linking number $\{\mu(rst);\, 1\leq
r<s<t\leq k\}$. 
The conjugations and partial conjugations act on this quotient by translations 
whose translation vectors' coordinates are linear functions of the linking 
numbers $l_{ij}$.
Thus, if the dimension of the subspace generated by these translation vectors
is less then $\binom{k}{3}$ for generic values of the linking numbers, we may
find a non-trivial vector perpendicular to all these translation vectors. Furthermore,
the coordinates of this vector could be taken as polynomials in $l_{ij}$.
Then the projection of an vector $\{\mu(rst)\}$ to this perpendicular 
vector will 
be invariant under conjugations and partial conjugations. This is the general 
philosophy behind our construction of link-homotopy invariant polynomials. 

In a recent preprint \cite{MT}, B. Mellor and D. Thurston has
established the existence 
of link-homotopy invariants of finite type which are not polynomials of
linking number when $k\geq9$. Their proof is not constructive and 
therefore it is
not clear whether their link-homotopy invariants are polynomials of 
$\mu$-numbers. 

On the other hand, since $\mu$-numbers are of finite type for string links 
\cite{L1,BN}, it is easy to see that our link-homotopy invariant
polynomials are of finite type for links. For $k\leq 5$, it is shown in 
\cite{MT} that the only finite type link-homotopy invariants are
polynomials in the linking numbers. So our construction fits nicely with
this work of Mellor and Thurston.

Recall that the only 
finite type knot concordance invariant is the Arf invariant \cite{Ng}. 
Since link
concordance implies link-homotopy, our work (as well as the work of Mellor 
and Thurston, of course) shows the existence of non-trivial finite
type link concordance invariants. 

To extend the applicability of our general philosophy slightly, we find that 
the
operation on the vector $\{\mu(rst)\}$ induced by reversing the orientation of 
each component of a string link is to change it by a negative sign followed
by a translation whose translation vector's coordinates are 
quadratic polynomials in $l_{ij}$.
If the dimension of the subspace generated by this vector together with
the translation vectors of conjugations and partial conjugations is still
less then $\binom{k}{3}$ for generic values of the linking numbers, 
and this is the case indeed, we can construct
a non-trivial link-homotopy invariant polynomial which is changed by a 
sign when
the orientation of each component of a link is reversed. 
We say that such a link invariant detects the invertibility for links. 
Recall that the reversion of the orientation of every component 
of a link does not change the quantum invariant associated with an irreducible
representation of a semi-simple Lie algebra (see, for example, \cite{L2}).
Thus our invariant is of finite type but is not determined by 
quantum invariants. The existence of a finite type knot invariant which
detects the invertibility for knots is a major problem in 
the theory of finite type invariants (see, for example, \cite{L2,K}).
We believe that finite type knot invariants can not detect the invertibility
for knots.

It remains unclear whether we can have a complete set of link-homotopy 
invariant polynomials which determines uniquely link-homotopy classes of 
links. See \cite{Le} for an earlier attempt on this problem.  
This problem could probably be translated to the problem of understanding 
the sublattice
generated by the translation vectors of conjugations and partial conjugations. 
A better understanding of this sublattice might also be useful 
in answering the following
question. If we let $\text{deg}\,(l_{ij})=1$ and $\text{deg}\,(\mu(rst))=2$, 
the link-homotopy invariant polynomial for $k=6$ we construct in Section 3, 
which detects the invertibility for links, is a linear combination of 
113,700 monomials of degree 22, homogeneous in both $l_{ij}$ and
$\mu(rst)$ and linear in $\mu(rst)$. Is there a 
shorter link-homotopy invariant polynomial detecting the invertibility for 
links?

\section{Conjugation and Partial Conjugation}

We first recall the classification of ordered, oriented links up to link-homotopy 
given in \cite{HL}. We will follow the notations of \cite{HL}.

Let $\H(k)$ be the group of link-homotopy classes of ordered, oriented
string links with
$k$ components. The components of a string link will be ordered by
$1,2,\dots,k$. It is shown in \cite{HL} that the pure braid group $P(k)$
of $k$ components maps onto $\H(k)$ under the natural map $P(k)\rightarrow
\H(k)$.

Deletion of the $i$-th component of the string link gives rise
to a group homomorphism $d_i:\H(k)\rightarrow\H(k-1)$. If $F(k)$ denotes the free
group of rank $k$ generated by $x_1,x_2,\dots,x_k$, the reduced free group
$RF(k)$ is the quotient of $F(k)$ by adding relations $[x_i,x_i^g]=1$ for
all $i$ and all $g\in F(k)$. 

\begin{lm} There is a split short exact sequence of groups
\begin{equation}\label{split}
1\rightarrow RF(k-1)\longrightarrow\H(k)\overset{d_i}
\longrightarrow\H(k-1)\rightarrow1
\end{equation}
where $RF(k-1)$ is the reduced free group generated by $x_1,\dots,x_{i-1},
x_{i+1},\dots,x_k$.
\end{lm}

Notice that the split exact sequence (\ref{split}) depends on the deleting 
component so that there are $k$ such split exact sequences all together. A split 
exact sequence determines a semi-direct product decomposition
$$\H(k)=\H(k-1)\ltimes RF(k-1).$$

Conjugation in the group $\H(k)$ is defined as usual: A conjugation of
$\sigma\in\H(k)$ by $\beta\in\H(k)$ is the element $\beta\sigma\beta^{-1}
\in\H(k)$. A {\it partial conjugation} of $\sigma\in\H(k)$ is 
an element of the form $\theta hgh^{-1}$, 
where we write $\sigma=\theta g$ according to a decomposition $\H(k)=
\H(k-1)\ltimes RF(k-1)$, for $\theta\in\H(k-1)$ and $g\in RF(k-1)$, 
and for an arbitrary $h\in RF(k-1)$. 

It is easy to see that every link-homotopy class of ordered, oriented links 
with $k$ components
can be realized as the closure of an element in $\H(k)$. One of the
main results of \cite{HL} is the following classification theorem.

\begin{thm}\label{classification} Let $\sigma,\sigma'\in\H(k)$. Then the
closures of $\sigma$ and $\sigma'$ are link-homotopic as ordered, oriented
links if and only if there is a sequence $\sigma=\sigma_0,\sigma_1,\dots,
\sigma_n=\sigma'$ of elements of $\H(k)$ such that $\sigma_{j+1}$
is either a conjugation or a partial conjugation of $\sigma_j$.
\end{thm}
 
For a group $G$, we will denote by $G_n$ the $n$-th term of the
lower central series of $G$, i.e. $G_1=G$ and $G_{n+1}=[G_n,G]$, 
the normal subgroup of $G$ generated by elements of the form 
$[g,h]=ghg^{-1}
h^{-1}$ for all $g\in G_n$ and $h\in G$. 
A group 
$G$ is nilpotent of class $n$ if $G_{n+1}=1$ but $G_{n}\neq1$. We summarize
some known facts about the group structures of $\H(k)$ in the 
following lemma. 

\begin{lm} 1) $\H(k)$ is torsion free and nilpotent of class $k-1$;

2) Corresponding to a decomposition $\H(k)=\H(k-1)\ltimes RF(k-1)$, we have
$$H(k)_n=\H(k-1)_n\ltimes RF(k-1)_n.$$

3) $\H(k)_{n-1}/\H(k)_{n}$ is a free abelian group of rank 
$(n-2)!\binom{k}{n}$.
\end{lm}

For $\sigma\in\H(k)$, its image in $\H(k)/\H(k)_3$ can be described
by $\binom{k}{2}+\binom{k}{3}$ integers. These integers are linking numbers
$l_{ij}$, for $1\leq i<j\leq k$, and Milnor's triple linking numbers 
$\mu(rst)$, for 
$1\leq r<s<t\leq k$. We want to have them defined precisely 
and understand how they change when $\sigma$ is changed by a 
conjugation or a partial conjugation.

We will denote by $\tau_{rs}=\tau_{sr}$, for $1\leq r<s\leq k$, the pure 
braid depicted
in Figure 1. Let $\sigma\in\H(k)/\H(k)_3$. For $1\leq r<s<t\leq k$, 
after deleting all
components other than the $r,s,t$-th components, $\sigma$ can be written
in the following normal form 
\begin{equation}\label{nform}
\sigma=\tau_{rs}^\alpha\tau_{rt}^\beta\tau_{st}^\gamma[\tau_{rt},\tau_{st}]^
\delta,
\end{equation}
where $\alpha=l_{rs}$, $\beta=l_{rt}$, $\gamma=l_{st}$.
By definition, we have $\delta=\mu(rst)$ for $\sigma\in\H(k)$.

\bigskip
\centerline{\epsfxsize=2.5in\epsfbox{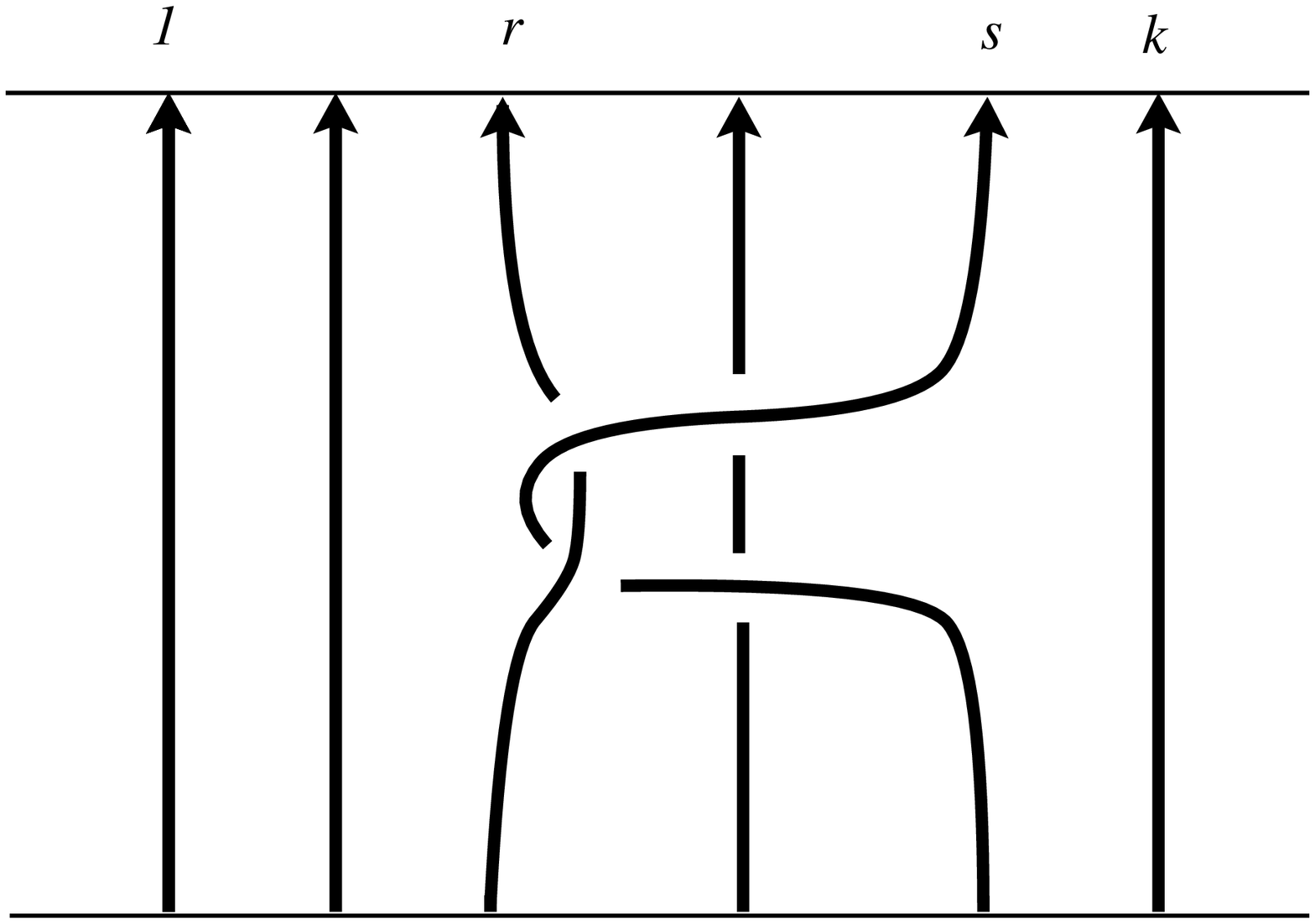}}
\medskip
\centerline{\small Figure 1. The pure braid $\tau_{rs}$.}
\bigskip
 
\begin{lm} In $\H(k)/\H(k)_3$, if $r',s',t'$ is a permutation of $r,s,t$ and
$\epsilon$ is the sign of the permutation, then
$$[\tau_{r't'},\tau_{s't'}]=[\tau_{rt},\tau_{st}]^{\epsilon}.$$
Furthermore, we have 
$$[\tau_{rt}^\eta,\tau_{st}]=[\tau_{rt},\tau_{st}]^\eta.$$
\end{lm}

This lemma is useful in the following calculation and its proof is straightforward.

To understand how $\mu(rst)$ changes under the conjugation, we only
need to calculate the conjugation of $\sigma\in\H(k)/\H(k)_3$ 
under the normal form (\ref{nform}) by
$\tau_{rs},\tau_{rt},\tau_{st}$. This calculation is straightforward:
$$
\begin{aligned}
\,&\begin{aligned}
\tau_{rs}\sigma\tau_{rs}^{-1}&=\tau_{rs}\tau_{rs}^\alpha\tau_{rt}^\beta
\tau_{st}^\gamma[\tau_{rt},\tau_{st}]^\delta\tau_{rs}^{-1}\\
&=\tau_{rs}^\alpha\tau_{rt}^\beta\tau_{st}^\gamma[\tau_{rs},\tau_{rt}]^{\beta}
[\tau_{rs},\tau_{st}]^{\gamma}[\tau_{rt},\tau_{st}]^\delta\\
&=\tau_{rs}^\alpha\tau_{rt}^\beta\tau_{st}^\gamma[\tau_{rt},\tau_{st}]^
{\delta+\beta-\gamma};
\end{aligned}\\
&\begin{aligned}
\tau_{rt}\sigma\tau_{rt}^{-1}&=\tau_{rt}\tau_{rs}^\alpha\tau_{rt}^\beta
\tau_{st}^\gamma[\tau_{rt},\tau_{st}]^\delta\tau_{rt}^{-1}\\
&=\tau_{rs}^\alpha\tau_{rt}^\beta\tau_{st}^\gamma[\tau_{rt},\tau_{rs}]^{\alpha}
[\tau_{rt},\tau_{st}]^{\gamma}[\tau_{rt},\tau_{st}]^\delta\\
&=\tau_{rs}^\alpha\tau_{rt}^\beta\tau_{st}^\gamma[\tau_{rt},\tau_{st}]^
{\delta-\alpha+\gamma};
\end{aligned}\\
&\begin{aligned}
\tau_{st}\sigma\tau_{st}^{-1}&=\tau_{st}\tau_{rs}^\alpha\tau_{rt}^\beta
\tau_{st}^\gamma[\tau_{rt},\tau_{st}]^\delta\tau_{st}^{-1}\\
&=\tau_{rs}^\alpha\tau_{rt}^\beta\tau_{st}^\gamma[\tau_{st},\tau_{rs}]^{\alpha}
[\tau_{st},\tau_{rt}]^{\beta}[\tau_{rt},\tau_{st}]^\delta\\
&=\tau_{rs}^\alpha\tau_{rt}^\beta\tau_{st}^\gamma[\tau_{rt},\tau_{st}]^
{\delta+\alpha-\beta}.
\end{aligned}
\end{aligned}$$

We summarize the calculation into the following lemma.

\begin{lm}\label{conj} The change of $\mu(rst)$ under a conjugation is given by
$$\begin{aligned}
&\text{Conjugation by $\tau_{rs}$:}\quad\mu(rst)\rightarrow
\mu(rst)+l_{rt}-l_{st};\\
&\text{Conjugation by $\tau_{rt}$:}\quad\mu(rst)\rightarrow
\mu(rst)-l_{rs}+l_{st};\\
&\text{Conjugation by $\tau_{st}$:}\quad\mu(rst)\rightarrow
\mu(rst)+l_{rs}-l_{rt}.
\end{aligned}
$$
Furthermore, $\mu(rst)$ will not change under a 
conjugation by $\tau_{ij}$ where
$\{i,j\}$ and $\{r,s,t\}$ have at most one element in common.
\end{lm}
 
The calculation of partial conjugations is slightly more complicated. We
will start with partial conjugations by $\tau_{rt}$ and $\tau_{st}$. 
These two operations are denoted by $\mathbf{t}^r$ and $\mathbf{t}^s$,
respectively. For $\sigma\in\H(k)/\H(k)_3$ under the normal form 
(\ref{nform}), we have:
$$
\begin{aligned}
\,&\begin{aligned}
\sigma\overset{\mathbf{t}^r}\longrightarrow&\,\,\tau_{rs}^\alpha\tau_{rt}
\tau_{rt}^\beta\tau_{st}^\gamma[\tau_{rt},\tau_{st}]^\delta\tau_{rt}^{-1}\\
&\,\,=\tau_{rs}\alpha\tau_{rt}^\beta\tau_{st}^\gamma[\tau_{rt},\tau_{st}]^{\delta+\gamma};
\end{aligned}\\
&\begin{aligned}
\sigma\overset{\mathbf{t}^s}\longrightarrow&\,\,\tau_{rs}^\alpha\tau_{st}
\tau_{rt}^\beta\tau_{st}^\gamma[\tau_{rt},\tau_{st}]^\delta\tau_{st}^{-1}\\
&\,\,=\tau_{rs}\alpha\tau_{rt}^\beta\tau_{st}^\gamma[\tau_{rt},\tau_{st}]^{\delta-\beta}.
\end{aligned}
\end{aligned}
$$

To calculate partial conjugations by $\tau_{rs}$ and $\tau_{ts}$, which 
are denoted by $\mathbf{s}^r$ and $\mathbf{s}^t$, respectively, we need to 
rewrite $\sigma$ as follows:
$$\sigma=\tau_{rs}^\alpha\tau_{rt}^\beta\tau_{st}^\gamma[\tau_{rt},\tau_{st}]
^\delta=\tau_{rt}^\beta\tau_{rs}^\alpha\tau_{ts}^\gamma[\tau_{rs},\tau_{ts}]
^{-\delta-\alpha\beta}.$$
Then, we have:
$$
\begin{aligned}
\,&\begin{aligned}
\sigma\overset{\mathbf{s}^r}\longrightarrow&\,\,\tau_{rt}^\beta\tau_{rs}
\tau_{rs}^\alpha\tau_{ts}^\gamma[\tau_{rs},\tau_{ts}]^{-\delta-\alpha\beta}
\tau_{rs}^{-1}\\
&\,\,=\tau_{rt}^\beta\tau_{rs}^\alpha\tau_{ts}^\gamma[\tau_{rs},\tau_{ts}]
^{-\delta-\alpha\beta+\gamma}\\
&\,\,=\tau_{rs}^\alpha\tau_{rt}^\beta\tau_{st}^\gamma[\tau_{rt},\tau_{st}]
^{\delta-\gamma};
\end{aligned}\\
&\begin{aligned}
\sigma\overset{\mathbf{s}^t}\longrightarrow&\,\,\tau_{rt}^\beta\tau_{ts}
\tau_{rs}^\alpha\tau_{ts}^\gamma[\tau_{rs},\tau_{ts}]^{-\delta-\alpha\beta}
\tau_{ts}^{-1}\\
&\,\,=\tau_{rt}^\beta\tau_{rs}^\alpha\tau_{ts}^\gamma[\tau_{rs},\tau_{ts}]
^{-\delta-\alpha\beta-\alpha}\\
&\,\,=\tau_{rs}^\alpha\tau_{rt}^\beta\tau_{st}^\gamma[\tau_{rt},\tau_{st}]
^{\delta+\gamma}.
\end{aligned}
\end{aligned}
$$

Similarly, to calculate partial conjugations $\mathbf{r}^s$ and $\mathbf{r}^t$,
we first rewrite $\sigma$:
$$\sigma=\tau_{rs}^\alpha\tau_{rt}^\beta\tau_{st}^\gamma[\tau_{rt},\tau_{st}]
^\delta=\tau_{st}^\gamma\tau_{sr}^\alpha\tau_{tr}^\beta
[\tau_{sr},\tau_{tr}]^{\delta-\alpha\gamma+\beta\gamma}.$$
Then, we have
$$
\begin{aligned}
\,&\begin{aligned}
\sigma\overset{\mathbf{r}^s}\longrightarrow&\,\,\tau_{st}^\gamma\tau_{sr}
\tau_{sr}^\alpha\tau_{tr}^\beta[\tau_{sr},\tau_{tr}]^{\delta-\alpha\gamma+
\beta\gamma}
\tau_{sr}^{-1}\\
&\,\,=\tau_{st}^\gamma\tau_{sr}^\alpha\tau_{tr}^\beta[\tau_{sr},\tau_{tr}]
^{\delta-\alpha\gamma+\beta\gamma+\beta}\\
&\,\,=\tau_{rs}^\alpha\tau_{rt}^\beta\tau_{st}^\gamma[\tau_{rt},\tau_{st}]
^{\delta+\beta};
\end{aligned}\\
&\begin{aligned}
\sigma\overset{\mathbf{r}^t}\longrightarrow&\,\,\tau_{st}^\gamma\tau_{tr}
\tau_{sr}^\alpha\tau_{tr}^\beta[\tau_{sr},\tau_{tr}]^{\delta-\alpha\gamma+
\beta\gamma}
\tau_{tr}^{-1}\\
&\,\,=\tau_{st}^\gamma\tau_{sr}^\alpha\tau_{tr}^\beta[\tau_{sr},\tau_{tr}]
^{\delta-\alpha\gamma+\beta\gamma-\alpha}\\
&\,\,=\tau_{rs}^\alpha\tau_{rt}^\beta\tau_{st}^\gamma[\tau_{rt},\tau_{st}]
^{\delta-\alpha}.
\end{aligned}
\end{aligned}
$$

We summarize the previous calculation into the following lemma.

\begin{lm}\label{pconj} 
The change of $\mu(rst)$ under a partial conjugation is given by
$$\begin{aligned}
&\mathbf{t}^r:\quad\mu(rst)\rightarrow\mu(rst)+l_{st};\\
&\mathbf{t}^s:\quad\mu(rst)\rightarrow\mu(rst)-l_{rt};\\
&\mathbf{s}^r:\quad\mu(rst)\rightarrow\mu(rst)-l_{st};\\
&\mathbf{s}^t:\quad\mu(rst)\rightarrow\mu(rst)+l_{rs};\\
&\mathbf{r}^s:\quad\mu(rst)\rightarrow\mu(rst)+l_{rt};\\
&\mathbf{r}^t:\quad\mu(rst)\rightarrow\mu(rst)-l_{rs}.
\end{aligned}
$$
Furthermore, a partial conjugation by $\mathbf{i}^j$ will not change 
$\mu(rst)$
if $\{i,j\}$ and $\{r,s,t\}$ have at most one element in common.
\end{lm}

For a given string link $\sigma\in\H(k)$, we will think of the whole 
collection $\{\mu(rst)\,;\,1\leq r<s<t\leq k\}$ as an element
in $\Z^{\binom{k}{3}}$. Then the conjugations and partial conjugations 
act on $\Z^{\binom{k}{3}}$ by translations. We will abuse the notation by
using the same symbol to denote both a translation operation and the 
corresponding translation vector. Thus, a translation operation 
$T:V\rightarrow V$ on a vector space $V$ is given by 
$T(v)=v+T$, for all $v\in V$ and a fixed $T\in V$. If $T_1$ and $T_2$
are two translations, we have 
$$(T_1\cdot T_2)(v)=v+T_1+T_2,\qquad\text{for all $v\in V$.}$$

The following 
two theorems follow directly from Lemmas \ref{conj}
and \ref{pconj}. 

\begin{thm}\label{translation} 
The translation operation on $\Z^{\binom{k}{3}}$ given by the
conjugation of $\tau_{ij}$ is the same as the composition of
the translation operations given by the partial conjugations $\mathbf{i}^j$
and $\mathbf{j}^i$, i.e. it is equal to $\mathbf{i}^j+\mathbf{j}^i$.
\end{thm} 

\begin{thm}\label{equation} The translation operations $\mathbf{i}^j$ 
satisfy the following relations:
$$\begin{aligned}
&\sum_{j\neq i}\mathbf{j}^i=0\\
&\sum_{j\neq i}l_{ij}\,\mathbf{i}^j=0
\end{aligned}$$
for all $i=1,2,\dots,k$.
\end{thm}

String links are oriented in the sense that each component is given a 
orientation from the bottom to the top. See Figure 1.
Reversing the orientation on each
component of a string link defines bijection 
$$\sigma\mapsto\overline{\sigma}:\H(k)\rightarrow\H(k).$$
This bijection is an anti-homomorphism: $\overline{\sigma_1\sigma_2}=
\overline{\sigma_2}\,\overline{\sigma_1}$.
This bijection induces an operation on $\Z^{\binom{k}{3}}$.

\begin{thm}\label{reversion} The operation on $\Z^{\binom{k}{3}}$ 
induced by reversing the 
orientation of each component of a string link is to change each
$\mu(rst)$ to $-\mu(rst)$ followed by the translation operation
$$\mu(rst)\longrightarrow\mu(rst)-l_{rs}\,l_{rt}+
l_{rs}\,l_{st}-l_{rt}\,l_{st}.
$$
\end{thm}

\begin{proof} Consider the normal form (\ref{nform}) of
$\sigma\in\H(k)/\H(k)_3$ in the $r,s,t$-th components. The normal
form for $\overline{\sigma}$ is obtained as follows:
$$\begin{aligned}
\overline{\sigma}&=[\tau_{rt},\tau_{st}]^{-\delta}\tau_{st}^\gamma
\tau_{rt}^\beta
\tau_{rs}^\alpha\\
&=\tau_{rs}^\alpha\tau_{rt}^\beta\tau_{st}^\gamma
[\tau_{rt},\tau_{st}]^{-\delta-\alpha\beta+\alpha\gamma-\beta\gamma}.
\end{aligned}
$$
Thus the operation on $Z^{\binom{k}{3}}$ induced by 
$\sigma\mapsto\bar\sigma$ is given by
$$\mu(rst)\longrightarrow-\mu(rst)-l_{rs}\,l_{rt}+
l_{rs}\,l_{st}-l_{rt}\,l_{st}.$$
\end{proof}

\section{Construction of the Invariant}

By Theorems \ref{classification} and \ref{translation}, 
we shall look for polynomials in 
$l_{ij}$ and $\mu(rst)$ invariant
under the translation operations on $\{\mu(rst)\}\in\Z^{\binom{k}{3}}$
induced by partial conjugations. There are $k(k-1)$ partial conjugations
all together and their induced translations subject to $2k$ linear
equations given in Theorem \ref{equation}. If these equations are
linearly independent for generic values of $\{l_{ij}\}$, 
the sublattice of $Z^{\binom
{k}{3}}$ generated by the translation vectors of the partial conjugations
will be of dimension no larger than $k(k-1)-2k=k^2-3k$. 

\begin{lm} For $k>3$, the $2k$ equations in Theorem \ref{equation} are
linearly independent for generic values of $\{l_{ij}\}$.
\end{lm}

\begin{proof} We write the two sets of equations in Theorem \ref{equation} 
as follows:
$$\mathbf{1}^i+\mathbf{2}^i+\cdots+\mathbf{j}^i+\cdots+\mathbf{k}^i=0,\quad
j\neq i;$$ 
$$l_{i1}\mathbf{i}^1+l_{i2}\mathbf{i}^2+\cdots+l_{ij}\mathbf{i}^j+\cdots
+l_{ik}\mathbf{i}^k=0,\quad j\neq i,$$
for each $i=1,2,\dots,k$.

For generic values of $\{l_{ij}\}$, 
use the first $k-1$ equations the first set of $k$ equations, we can
solve for $\mathbf{k}^1,\mathbf{k}^2,\dots,\mathbf{k}^{k-1}$. Similarly,
we can solve for $\mathbf{1}^k,\mathbf{2}^k,\dots,(\mathbf{k-1})^k$ from
the first $k-1$ equations of the second set of $k$ equations. The remaining
vectors $\mathbf{i}^j$, $i,j\neq k$, have to satisfy another
two equations obtained from the last equations in those two sets of
$k$ equations, respectively, by substituting $\mathbf{k}^i$ and $\mathbf{i}^k$ with 
their solutions in terms of $\mathbf{i}^j$ for $i,j\neq k$. It it then easy 
to check that these two equations are linearly independent when $k>3$.
\end{proof}

\begin{lm} For $k=4,5$, we have $\binom{k}{3}=k^2-3k$. For $k\geq 6$,
we have $\binom{k}{3}>k^2-3k$.
\end{lm} 

\begin{proof} The case of $k=4,5$ can be checked directly. For $k\geq6$, we
have
$$\begin{aligned}
\binom{k}{3}-(k^2-3k)&=\frac{k}{6}\,(k^2-9k+20)\\
&=\frac{k}{6}\,((k-6)^2+3k-16)>0.
\end{aligned}
$$
\end{proof}

\begin{thm}\label{exist} For $k\geq 6$, there exists a polynomial in 
$l_{ij}$
and $\mu(rst)$ which is a link-homotopy invariant of ordered, oriented
links with $k$ components. This link-homotopy invariant is of finite type.
\end{thm}

\begin{proof} In $\Z^{\binom{k}{3}}$, let $\P$ be the sublattice generated
by the translation vectors of partial conjugations. Then we have
$$\text{dim}(\P)\leq k^2-3k<\binom{k}{3}.$$

Let $\Omega\in\Z^{\binom{k}{3}}$ be a non-zero vector perpendicular to 
$\P$. We can choose such a $\Omega$ so that its coordinates  
are polynomials in $\{l_{ij}\}$ and the inner product 
$\mathbf{i}^j\cdot\Omega$ is identically zero. This can be achieved by
considering generic values of $\{l_{ij}\}$ first. Then since
$\mathbf{i}^j\cdot\Omega=0$ for generic values of $\{l_{ij}\}$,
it has to be zero identically. Let
$\mu=\{\mu(rst)\}\in\Z^{\binom{k}{3}}$. The inner product $\mu\cdot\Omega$
is invariant under the translations by vectors in $\P$. This is a
desired link-homotopy invariant of ordered, oriented links since
$$(\mu+\mathbf{i}^j)\cdot\Omega=\mu\cdot\Omega$$
for all $i,j=1,2,\dots,k$.

The fact that the invariant $\mu\cdot\Omega$ is of fine type is a direct
consequence of the fact that the linking numbers and the triple linking numbers
are all finite type invariants of string links (\cite{L1,BN}). If we 
have a singular link, we may put it into the form of the closure of a single
string link. Since polynomials of finite type invariants are still of finite 
type, $\mu\cdot\Omega$ vanishes on singular string links with sufficiently large
number of double points. This implies that it is a finite type link invariant.
\end{proof} 
 
We now consider in some details the case $k=6$. 

Let us order $\mu(rst)$, $1\leq r<s<t\leq 6$ in the lexicographic order. So
$$\begin{aligned}
\mu=(&\mu(123),\mu(124),\mu(125),\mu(126),\mu(134),\mu(135),\mu(136),
\mu(145),\mu(146),\mu(156),\\
&\mu(234),\mu(235),\mu(236),\mu(245),\mu(246),\mu(256),\mu(345),\mu(346),
\mu(356),\mu(456)).
\end{aligned}
$$
Then the vectors of the translation operations $\mathbf{1}^2$, 
$\mathbf{1}^3$, $\mathbf{1}^4$,
$\mathbf{1}^5$,
$\mathbf{1}^6$,
$\mathbf{2}^1$, $\mathbf{2}^3$, $\mathbf{2}^4$, $\mathbf{2}^5$, 
$\mathbf{2}^6$, $\mathbf{3}^1$, $\mathbf{3}^2$,
$\mathbf{3}^4$, $\mathbf{3}^5$, $\mathbf{3}^6$, $\mathbf{4}^1$,
$\mathbf{4}^2$, $\mathbf{4}^3$, $\mathbf{4}^5$, $\mathbf{4}^6$, 
$\mathbf{5}^1$, $\mathbf{5}^2$,
$\mathbf{5}^3$,$\mathbf{5}^4$, $\mathbf{5}^6$,
$\mathbf{6}^1$, $\mathbf{6}^2$, $\mathbf{6}^3$, $\mathbf{6}^4$, 
$\mathbf{6}^5$  are the row vectors of the following $30\times 20$ matrix
from the top to the bottom, respectively:

$$\left|\begin{smallmatrix}\setcounter{MaxMatrixCols}{20}
 l_{13}& l_{14}& l_{15}& l_{16}& 0& 0& 0& 0& 0& 0& 0& 0& 0& 0& 0& 0& 0& 0& 0& 0 \\
-l_{12}& 0& 0& 0& l_{14}& l_{15}& l_{16}& 0& 0& 0& 0& 0& 0& 0& 0& 0& 0& 0& 0& 0 \\
 0&-l_{12}& 0& 0&-l_{13}& 0& 0& l_{15}& l_{16}& 0& 0& 0& 0& 0& 0& 0& 0& 0& 0& 0 \\
 0& 0&-l_{12}& 0& 0&-l_{13}& 0&-l_{14}& 0& l_{16}& 0& 0& 0& 0& 0& 0& 0& 0& 0& 0 \\
 0& 0& 0&-l_{12}& 0& 0&-l_{13}& 0&-L_{14}&-l_{15}& 0& 0& 0& 0& 0& 0& 0& 0& 0& 0 \\
-l_{23}&-l_{24}&-l_{25}&-l_{26}& 0& 0& 0& 0& 0& 0& 0& 0& 0& 0& 0& 0& 0& 0& 0& 0 \\
 l_{12}& 0& 0& 0& 0& 0& 0& 0& 0& 0& l_{24}& l_{25}& l_{26}& 0& 0& 0& 0& 0& 0& 0 \\
 0& l_{12}& 0& 0& 0& 0& 0& 0& 0& 0&-l_{23}& 0& 0& l_{25}& l_{26}& 0& 0& 0& 0& 0 \\
 0& 0& l_{12}& 0& 0& 0& 0& 0& 0& 0& 0&-l_{23}& 0&-l_{24}& 0& l_{26}& 0& 0& 0& 0 \\
 0& 0& 0& l_{12}& 0& 0& 0& 0& 0& 0& 0& 0&-l_{23}& 0&-l_{24}&-l_{25}& 0& 0& 0& 0 \\
 l_{23}& 0& 0& 0&-l_{34}&-l_{35}&-l_{36}& 0& 0& 0& 0& 0& 0& 0& 0& 0& 0& 0& 0& 0 \\
-l_{13}& 0& 0& 0& 0& 0& 0& 0& 0& 0&-l_{34}&-l_{35}&-l_{36}& 0& 0& 0& 0& 0& 0& 0 \\
 0& 0& 0& 0& l_{13}& 0& 0& 0& 0& 0& l_{23}& 0& 0& 0& 0& 0& l_{35}& l_{36}& 0& 0 \\
 0& 0& 0& 0& 0& l_{13}& 0& 0& 0& 0& 0& l_{23}& 0& 0& 0& 0&-l_{34}& 0& l_{36}& 0 \\
 0& 0& 0& 0& 0& 0& l_{13}& 0& 0& 0& 0& 0& l_{23}& 0& 0& 0& 0&-l_{34}&-l_{35}& 0 \\
 0& l_{24}& 0& 0& l_{34}& 0& 0&-l_{45}&-l_{46}& 0& 0& 0& 0& 0& 0& 0& 0& 0& 0& 0 \\
 0&-l_{14}& 0& 0& 0& 0& 0& 0& 0& 0& l_{34}& 0& 0&-l_{45}&-l_{46}& 0& 0& 0& 0& 0 \\
 0& 0& 0& 0&-l_{14}& 0& 0& 0& 0& 0&-l_{24}& 0& 0& 0& 0& 0&-l_{45}&-l_{46}& 0& 0 \\
 0& 0& 0& 0& 0& 0& 0& l_{14}& 0& 0& 0& 0& 0& l_{24}& 0& 0& l_{34}& 0& 0& l_{46} \\
 0& 0& 0& 0& 0& 0& 0& 0& l_{14}& 0& 0& 0& 0& 0& l_{24}& 0& 0& l_{34}& 0&-l_{45} \\
 0& 0& l_{25}& 0& 0& l_{35}& 0& l_{45}& 0&-l_{56}& 0& 0& 0& 0& 0& 0& 0& 0& 0& 0 \\
 0& 0&-l_{15}& 0& 0& 0& 0& 0& 0& 0& 0& l_{35}& 0& l_{45}& 0&-l_{56}& 0& 0& 0& 0 \\
 0& 0& 0& 0& 0&-l_{15}& 0& 0& 0& 0& 0&-l_{25}& 0& 0& 0& 0& l_{45}& 0&-l_{56}& 0 \\
 0& 0& 0& 0& 0& 0& 0&-l_{15}& 0& 0& 0& 0& 0&-l_{25}& 0& 0&-l_{35}& 0& 0&-l_{56} \\
 0& 0& 0& 0& 0& 0& 0& 0& 0& l_{15}& 0& 0& 0& 0& 0& l_{25}& 0& 0& l_{35}& l_{45} \\
 0& 0& 0& l_{26}& 0& 0& l_{36}& 0& l_{46}& l_{56}& 0& 0& 0& 0& 0& 0& 0& 0& 0& 0 \\ 
 0& 0& 0&-l_{16}& 0& 0& 0& 0& 0& 0& 0& 0& l_{36}& 0& l_{46}& l_{56}& 0& 0& 0& 0 \\
 0& 0& 0& 0& 0& 0&-l_{16}& 0& 0& 0& 0& 0&-l_{26}& 0& 0& 0& 0& l_{46}& l_{56}& 0 \\
 0& 0& 0& 0& 0& 0& 0& 0&-l_{16}& 0& 0& 0& 0& 0&-l_{26}& 0& 0&-l_{36}& 0& l_{56} \\
 0& 0& 0& 0& 0& 0& 0& 0& 0&-l_{16}& 0& 0& 0& 0& 0& l_{26}& 0& 0&-l_{36}&-l_{46} \\
\end{smallmatrix}
\right|
$$
\medskip

We shall pick out the 18 rows of this matrix corresponding to the translation
operations of 
$\mathbf{1}^2$, $\mathbf{1}^3$, $\mathbf{1}^4$, $\mathbf{1}^5$,
$\mathbf{2}^1$, $\mathbf{2}^3$, $\mathbf{2}^4$, $\mathbf{2}^5$, 
$\mathbf{3}^1$, $\mathbf{3}^2$, $\mathbf{3}^4$, $\mathbf{3}^5$, $\mathbf{4}^1$,
$\mathbf{4}^2$, $\mathbf{4}^3$, $\mathbf{4}^5$, $\mathbf{5}^1$, $\mathbf{5}^2$,
respectively. Calculation using Mathematica shows that these 18 vectors are
linearly independent generically. 

Consider now the operation of reversing the orientation. The vector $R=\{R(rst)\}
\in\Z^{20}$
of the translation operation in Theorem \ref{reversion} is given by
$$R(rst)=-l_{rs}\,l_{rt}+l_{rs}\,l_{st}-l_{rt}\,l_{st}.$$
One can verify that the vector $R$ and the previous 18 vectors 
are linearly independent. Let $\M$ be the $19\times 20$ matrix formed by
these 19 vectors. Let $\M^{(i)}$ be the $19\times19$ matrix obtained from $\M$ by
deleting the $i$-th column from $\M$, $i=1,2,\dots,20$. Let
$$\Omega_i=(-1)^{i-1}\text{det}\,(\M^{(i)})$$
and $\Omega=(\Omega_1,\Omega_2,\dots,\Omega_{20})$.

\begin{thm} $\mu\cdot\Omega$ is a finite type link-homotopy invariant of ordered, 
oriented links with 6 component. When the orientation of every component is reversed,
this invariant is changed only by a sign. 
\end{thm}

\begin{proof} Using the fact that the rows of the cofactor matrix $A^*$ of a given 
matrix $A$ are perpendicular to different rows of $A$, we see that $\Omega$ is 
perpendicular to all the vectors of translation operation induced by partial 
conjugations as well as the vector $R$. Certainly, $\Omega\neq0$. So 
$\mu\cdot\Omega$ is a non-trivial link-homotopy 
invariant of ordered, oriented
links with 6 components. It is of finite type since it is a polynomial in $l_{ij}$
and $\mu(rst)$. Under the reversion of orientation, $\mu$ changes to $-\mu+R$.
Since $R\cdot\Omega=0$, the invariant $\mu\cdot\Omega$ is only changed by a sign under
the reversion of orientation. 
\end{proof}

To finish, let us furnish some data obtained using Mathematica. 
Let $\text{deg}\,(l_{ij})=1$, then $\Omega_i$ is a homogeneous polynomial
of degree 20 in $l_{ij}$. Let $L_i$ be the number of monomials in $\Omega_i$,
the sequence $\{L_1,L_2,\dots,L_{20}\}$ is given as follows:
$$\begin{aligned}
\{&5531, 5555, 5555, 5531, 5424, 5769, 5802, 5734, 5753, 5432,\\ 
&5432, 5753, 5802, 5734, 5769, 5424, 5928, 5922, 5922, 5928\}.
\end{aligned}
$$
Thus $\mu\cdot\Omega$ is linear and homogeneous in $\mu(rst)$ and has
113,700 monomials.

\end{document}